\documentclass{pnastwo}
\usepackage{t1enc} 

\usepackage[dvips]{graphicx}
\usepackage{pnastwoF}
\usepackage{amssymb,amsfonts,amsmath}
\usepackage{color}

\newcommand{\R}{{\mathbb R}}
\newcommand{\C}{{\mathbb C}}

 \newcommand{\Id}{{\mathbf I}}
 
\renewcommand{\O}{\mbox{\bf O}} 
\newcommand{\SO}{\mbox{\bf SO}}

\newcommand{\E}{\mbox{\bf E}}

\copyrightyear{2008}
\issuedate{Issue Date}
\volume{Volume}
\issuenumber{Issue Number}

\begin{document}

\title{A Huygens principle for diffusion and anomalous diffusion in spatially extended systems}

\author{Georg A. Gottwald\affil{1}{School of Mathematics and Statistics, University of Sydney, Sydney 2006 NSW, Australia}
\and
Ian Melbourne\affil{2}{Department of Mathematics, University of Surrey, Guildford GU2 7XH, U.K. (Current address: Mathematics Institute, University of Warwick, Coventry CV4 7AL, U.K.)}}

\contributor{Submitted to Proceedings of the National Academy of Sciences of the United States of America}

\maketitle

\begin{article}
\begin{abstract}
We present a universal view on diffusive behaviour in chaotic spatially extended systems for anisotropic and isotropic media. For anisotropic systems, strong chaos leads to diffusive behaviour (Brownian motion with drift) and weak chaos leads to superdiffusive behaviour (L\'evy processes with drift). For isotropic systems, the drift term vanishes and strong chaos again leads to Brownian motion. We establish the existence of a nonlinear Huygens principle for weakly chaotic systems in isotropic media whereby the dynamics behaves diffusively in even space dimension and exhibits superdiffusive behaviour in odd space dimensions.
\end{abstract}

\keywords{dynamical systems | anomalous diffusion | pattern formation | lattice models | symmetry}



\dropcap{I}n this article we study the dynamics of spatially extended systems with symmetry. In particular we consider anisotropic systems which are invariant under translations but not rotations, and isotropic systems which are invariant under translations and rotations. Examples of anisotropic systems are those with a preferred spatial direction such as fluid systems advected by a directed mean flow or reaction-advection diffusion systems; examples of isotropic systems are those without any preferred direction, e.g. reaction-diffusion systems. This paper is concerned with the following general question: Given a certain type of dynamics (regular or chaotic), can we say anything about the diffusive behaviour of the solution.

The answer to this question is different for isotropic and nonisotropic systems and is summarised in Table~1.   As we explain below, many of these results exist in some form in the literature and our main contribution there is to bring them together in a unified way. Furthermore, we establish a dichotomy for weakly chaotic dynamics in isotropic systems, whereby generically there is diffusion in even dimensions but superdiffusion in odd dimensions, which does not seem to have been noticed previously.

In the remainder of this section, we describe the results summarised in Table~1. Throughout, we let $d$ denote the number of space dimensions. In anisotropic systems where there is translation symmetry but no rotation symmetry, the simplest kind of solution is a travelling wave propagating with linear speed $c\in\R^d$.   By regular dynamics we mean travelling waves and periodic modulations of travelling waves; such solutions are bounded in a frame of reference moving with constant speed. Coullet and Emilsson \cite{CoulletEmilsson96} considered a family of modified Ginzburg-Landau equations on the line ($d=1$) and observed that chaotically modulated travelling waves exhibit Brownian motion-like diffusive behaviour in a frame of reference moving with constant wave speed.  The situation can be generalized to arbitrary dimensions and is made rigorous in~\cite{MelbourneNicol09} for strongly chaotic systems. The term {\em strong chaos} is defined below and includes many of the classical examples of chaotic systems. For example, systems satisfying the Gallavotti-Cohen Chaotic Hypothesis are strongly chaotic, as is the classical Lorenz attractor. Recently interest has grown in so called {\em{weakly chaotic}} dynamics, exemplified by Pomeau-Manneville intermittency~\cite{PomeauManneville80}. In the anisotropic case, this leads to superdiffusive behaviour~\cite{GaspardWang88,Gouezel04,MelbourneZweimueller13}. 

The isotropic case, where the dynamical system is invariant under the Euclidean group of rotations and translations of $d$-dimensional space, $d\ge2$, exhibits a richer variety of possible behaviours. The additional invariance with respect to rotations may lead to cancellations of the linear drift present in the anisotropic cases. For regular dynamics, a dichotomy between even and odd spatial dimensions exists. Analogously to the Huygens principle which states that one can hear only in an ambient odd dimensional space, here a nonlinear Huygens principle is operating: in even dimensions the solutions are bounded, whereas in odd dimensions solutions propagate linearly. A manifestation of this is the behaviour of spiral waves in $2$-dimensional excitable media \cite{Winfree91,Barkley94}; the spiral tip moves around in circles or flower-petal meanders. In the strongly chaotic case, the dichotomy disappears and the behaviour is like a Brownian motion without drift. (In the situation of spiral waves, such behaviour is called {\em hypermeander}~\cite{Rossler79,Winfree91,BiktashevHolden98,AshwinMelbourneNicol01}.   Although there is now a good theoretical understanding~\cite{NicolMelbourneAshwin01,FMT03,MelbourneTorok04,MelbourneNicol04}, a conclusive demonstration of the existence of hypermeander in physical or numerical investigations of spiral wave dynamics remains an open problem.) The weakly chaotic case in the isotropic case has not been previously studied.  Again there is no linear drift.  We will establish here another instance of a nonlinear Huygens principle: generically superdiffusive behaviour prevails in odd dimensions but the superdiffusion is suppressed in even dimensions and replaced by Brownian motion.

\begin{table}[h]
 \begin{tabular}{@{\vrule height 10.5pt depth4pt  width0pt}clll} 
\hline
Dynamics & Anisotropic medium & \multicolumn{2}{c}{Isotropic medium}\\
& \multicolumn{1}{c}{$d\ge1$} & \multicolumn{1}{c}{$d\ge3$ odd} & $d\ge2$ even
\\ 
\hline
Regular & $ct+\text{bounded}$ & $ct+\text{bounded}$ & bounded
\\ 
\hline
Strongly chaotic & $ct+\text{diffusive}$
& diffusive & diffusive 
\\ 
\hline
Weakly chaotic & $ct+\text{superdiffusive}$  &
superdiffusive  & diffusive
\\ 
\hline
\end{tabular}
\label{table}

\end{table}
\noindent\textbf{Table 1.} Rates of propagation for given dynamics in a $d$-dimensional anisotropic or isotropic medium. $c\in\R^d$ denotes a general vector in $\R^d$


\section{Spatially extended systems with symmetry}
We adopt here the standard perspective of decomposing the dynamics into the dynamics on the symmetry group and the dynamics orthogonal to it. Systems with symmetry or {\em{equivariant dynamical systems}} are thus cast into a so called {\em{skew product}} of the form 
\begin{align} \label{eq-skew}
\dot x=f(x), \quad  \dot g=g\xi(x),
\end{align}
on $X\times G$, where the dynamics on the symmetry group $G$ is driven by the shape dynamics on a cross-section $X$ transverse to the group directions. Here, $g\xi(x)$ denotes the action of the group element $g\in G$ on $\xi(x)\in T_eG$ (the Lie algebra of $G$).
Substituting the solution $x(t)$ for the shape dynamics into the $\dot g$ equation yields the nonautonomous equation $\dot g=g\xi(x(t)) $ to be solved for the group dynamics.  
The simplest example is the case of a travelling wave in a one-dimensional system with translation symmetry, where the shape dynamics is an equilibrium solution $x(t)\equiv x_0$ in the frame of reference moving with constant wave speed, and the dynamics on the translation group orbit describes the linear drift of the reference frame in physical space. 
A more interesting example is provided by spiral waves in excitable media where periodic shape dynamics leads to {\em meandering} of the spiral tip \cite{Winfree72}.  Here the shape dynamics $x(t)$ is periodic and the group dynamics $g(t)$ evolves quasiperiodically.

The skew product formulation has proved 
successful in describing both local bifurcations and global dynamics in pattern formation \cite{Barkley94,BiktashevHolden98,AshwinMelbourneNicol01,FiedlerSandstedeScheelWulff96,SandstedeScheelWulff97,GolubitskyLeBlancMelbourne00}, constructing efficient numerical methods 
for equivariant systems \cite{BeynThuemmler04,HermannGottwald10,FoulkesBiktashev10}, and in studying Hamiltonian systems to explain for example dynamics of periodic orbits relevant to planetary dynamics \cite{RobertsWulffLamb02} and observed spectra of $CO_2$ molecules \cite{CushmanEtAl04}.

Adopting the decomposition~\eqref{eq-skew}, we can rephrase our main question: Given a certain shape dynamics for $\dot x=f(x)$ (regular, strongly chaotic or weakly chaotic), can we say anything about the growth rate of solutions $g(t)$ on the group. More specifically what is the expected diffusive behaviour of $g(t)$? In the following we first recall results for the anisotropic case, and then treat isotropic systems.


\subsection{Anisotropic case}
In the anisotropic case, the symmetry group is the group of translations $G=\R^d$ of $d$-dimensional space. The skew product~\eqref{eq-skew} reduces to
\begin{align}
\dot x &= f(x)
\label{e.aniso_x}
\\
\dot p &= \phi(x)\; ,
\label{e.aniso_p}
\end{align}
where $\phi$ takes values in $\R^d$ and $p\in\R^d$ represents the translation variable. Without loss of generality we assume initial conditions $x(0)=x_0$ and $p(0)=0$. Equation \eqref{e.aniso_p} can be integrated to yield
\begin{align}
p(t) = \int_0^t\phi(x(s))\, ds\; .
\label{e.p}
\end{align}
If the shape dynamics \eqref{e.aniso_x} consists of an equilibrium $x(t)\equiv x_0$ we obtain $p(t)=c t$ with $c=\phi(x_0)$. This includes the case of a travelling wave moving with constant speed $c$ mentioned above. For a periodic solution $x(t+T)=x(t)$, we obtain $p(t)=ct+O(1)$ with $c = \int_0^T\phi(x(s))\, ds$. 

Next suppose that there is a chaotic attractor $\Lambda\subset X$ for the shape dynamics with ergodic invariant measure $\mu$. Using \eqref{e.p} it follows from the Birkhoff ergodic theorem that for typical initial conditions $x_0$,
\begin{align*}
\frac{1}{t}p(t) = \frac{1}{t}\int_0^t\phi(x(s))\, ds \longrightarrow c
\end{align*}
where $c= \int_\Lambda \phi\, d\mu$ is the time-average of $\phi$. Typically $c\neq0$, in which case there is linear drift as for the regular case. For strongly chaotic shape dynamics, it follows from~\cite{MelbourneNicol05,MelbourneNicol09,Gouezel10} that there exists $\lambda>0$ such that for typical initial conditions $x_0$,
\begin{align} \label{eq-ASIP}
p(t) = ct+ W(t) + O(t^{\frac12-\lambda}) \qquad {\rm{a.e.}}
\end{align}
where $W$ is a $d$-dimensional Brownian motion with covariance matrix $\Sigma$. This implies the central limit theorem: $\mu(x_0:(p(t)-ct)/\sqrt{t}\in I)\to {\rm Pr}(Y\in I)$ for each rectangle $I\subset\R^d$, where $Y\sim N(0,\Sigma)$ is a normally distributed $d$-dimensional random variable with mean $0$ and covariance matrix $\Sigma$. Another consequence is that the sequence $(p(nt)-cnt)/\sqrt{n}$ converges weakly to $W$ in the space of continuous sample paths. (This is called weak convergence of $p(t)-ct$ to Brownian motion, whereas~\eqref{eq-ASIP} is strong convergence.) For $d=1$, this describes the case of chaotically modulated travelling waves as observed in \cite{CoulletEmilsson96}.

Weakly chaotic dynamical systems are characterized by ``sticky'' equilibria,  periodic solutions and so on, where the dynamics exhibits laminar behaviour interspersed with intermittent chaotic bursts~\cite{PomeauManneville80}. It is well-known~\cite{GaspardWang88} that for such intermittent systems the usual central limit theorem may break down leading to fluctuations of L\'evy type rather than of Gaussian type. In those situations it was established~\cite{Gouezel04} that solutions propagate superdiffusively as $t^\gamma$ for some $\gamma\in(\frac12,1)$.  More precisely, $t^{-\gamma}(p(t)-ct)$ converges in distribution to an $\alpha$-stable law where $\alpha=1/\gamma$. Let $W_\alpha$ denote the corresponding L\'evy process (possessing increments that are independent, stationary, and with distributions proportional to this stable law). Then $p(t)-ct$ converges weakly to $W_\alpha(t)$ by~\cite{MelbourneZweimueller13}. This concludes the discussion of the anisotropic case in Table~1.


\subsection{Isotropic case}
In the isotropic case, the symmetry group is the Euclidean group $\E(d)=\SO(d)\ltimes\R^d$ consisting of rotations and translations of $d$-dimensional space, and the skew product equations are given by
\begin{align}
\label{e.iso}
\dot x = f(x), \quad
\dot A = A h(x), \quad
\dot p = A v(x)\; ,
\end{align}
where $A\in \SO(d)$ represents the rotation variables and $p\in \R^d$ represents the translation variables. Without loss of generality we choose as initial conditions $x(0)=x_0$, $A(0)=\Id$ and $p(0)=0$. Note that $\SO(d)$ consists of $d\times d$ orthogonal matrices with determinant $1$, and that $h$, being an element of the Lie algebra of $\SO(d)$, is a skew-symmetric matrix. (We suppose throughout that $d\ge2$, since otherwise we are in the anisotropic situation.)

If the shape dynamics consists of an equilibrium $x(t)\equiv x_0$, then the dynamics on the rotation group can be integrated to yield $A(t)=\exp(t h(x_0))$. We choose coordinates so that the skew-symmetric matrix $h(x_0)$ is diagonal with entries on the imaginary axis. In even dimensions $d=2q$, the diagonal entries are given by $\pm i\omega_1,\ldots,\pm i\omega_q$ and are typically nonzero. Using the identification $\R^d\cong\C^{q}$, we obtain $\dot p_j=e^{it\omega_j}v_j(x_0)$ and hence 
\begin{align}
\label{e.peven}
p_j(t)=(1/i\omega_j)e^{it\omega_j}v_j(x_0),\enspace j=1,\dots,q. 
\end{align}
It follows that $p(t)$ is bounded.

In odd dimensions $d=2q+1$, one of the diagonal entries of $h(x_0)$, without loss the first entry, is forced to vanish and with the identification $\R^d\cong\R\times\C^q$ we have $\dot p_1=v_1(x_0)$. Hence 
\begin{align} \label{e.podd}
p(t)=ct+O(1), \quad c=(v_1(x_0),0,\dots,0),
\end{align}
and there is typically a linear drift. For periodic solutions, analogous calculations~\cite{AshwinMelbourne97} lead similarly to bounded motion for $d$ even and  unbounded linear growth for $d$ odd. This constitutes the nonlinear Huygens principle for regular dynamics in isotropic media. This dichotomy can be visualized by looking at the effect of the rotations in even and odd dimensions. In even dimensions, all components of $v$ in \eqref{e.iso} are rotated; in odd dimensions, however, there is an axis of rotation and the corresponding component of $v$ is not subjected to the averaging effect of the rotation. In Figure~1 we show typical dynamics of the translation variables for $\E(2)$ and $\E(3)$ skew products with underlying regular dynamics. We see clearly bounded motion for $d=2$ and a corkscrew motion along the axis of rotation of $\SO(3)$ for $d=3$.

\begin{figure}
\begin{center}
\centerline{
\includegraphics[width=.2\textwidth]{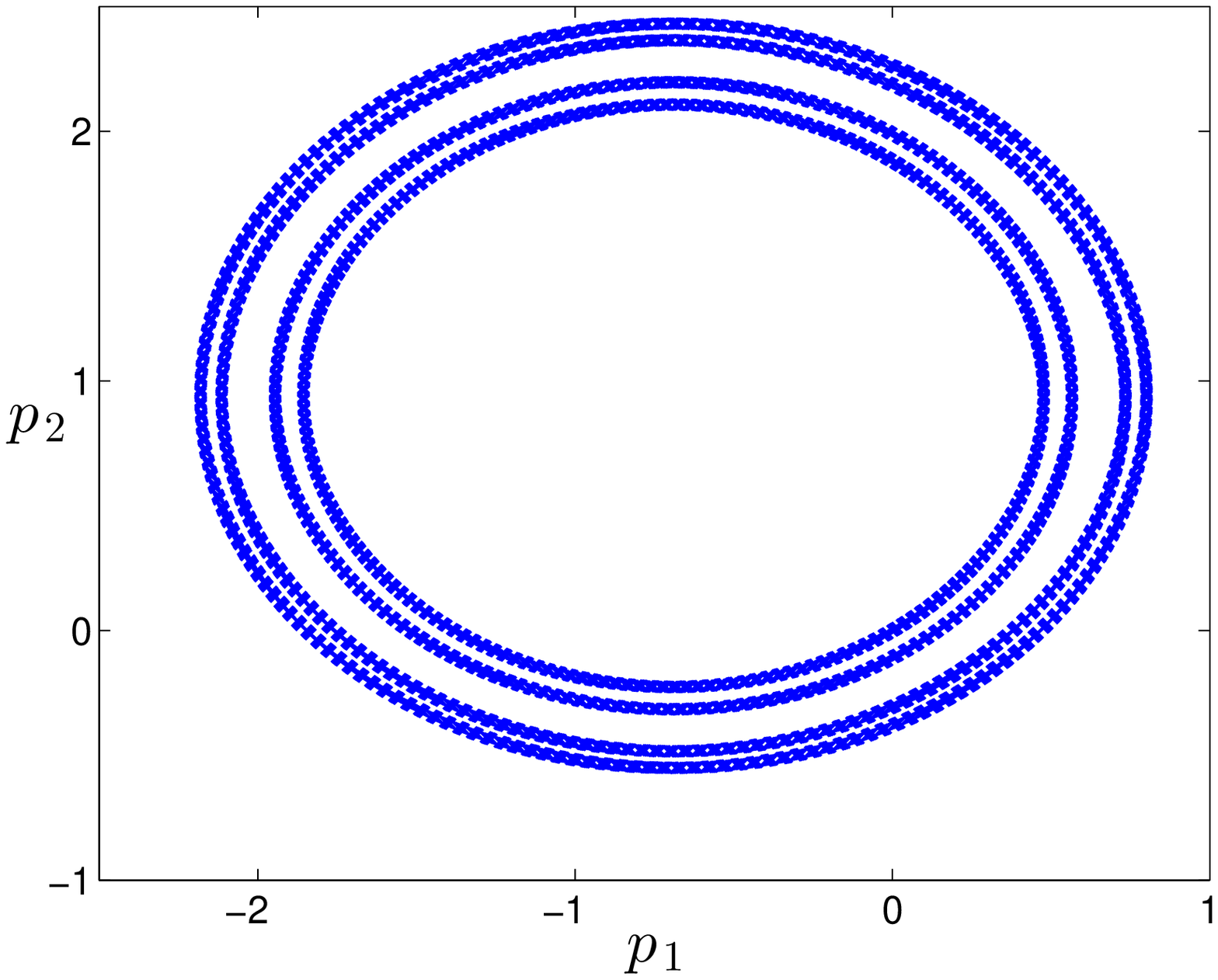}
\includegraphics[width=.2\textwidth]{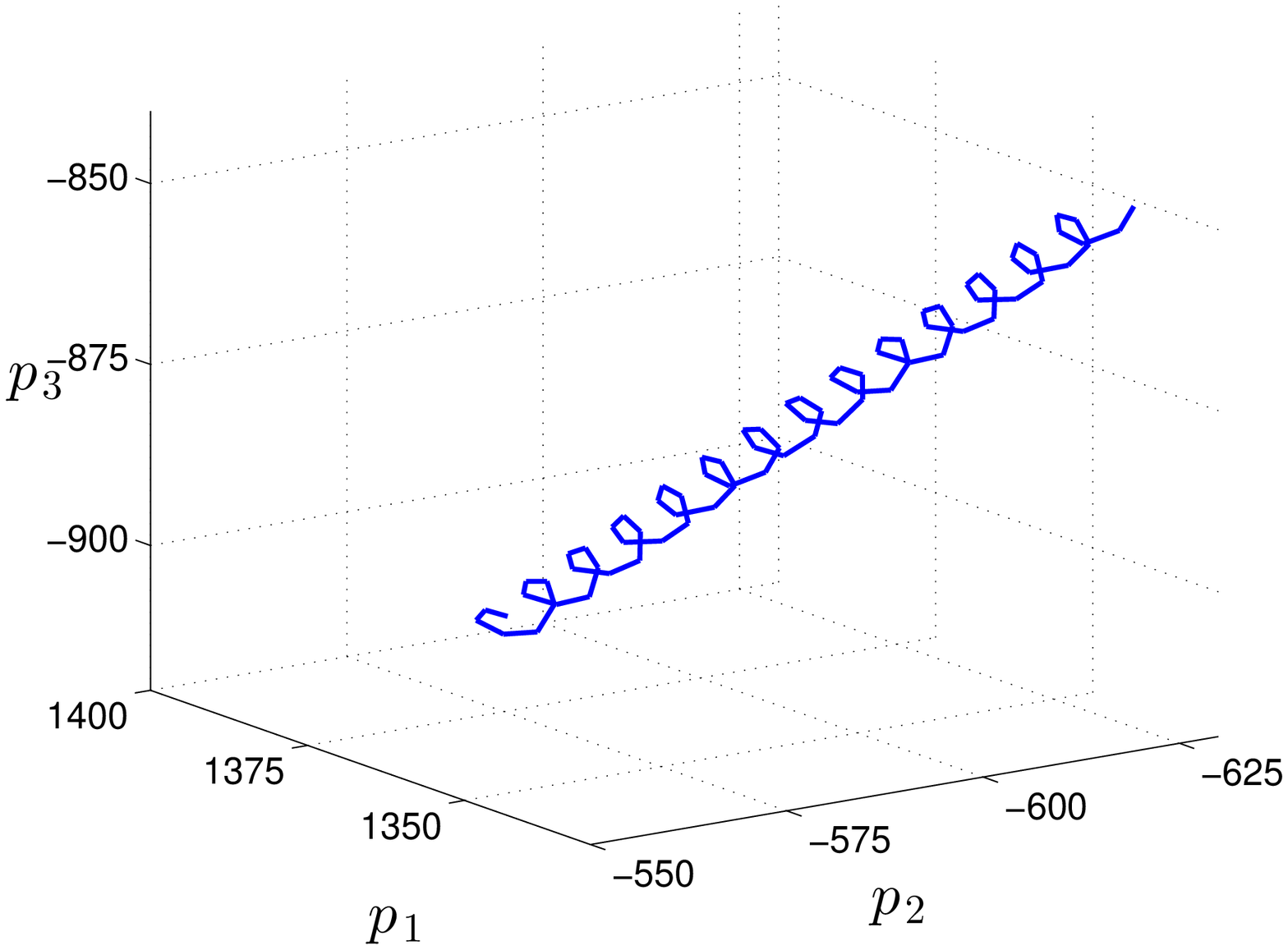}}
\caption{Typical behaviour of the translation variables $p$ for underlying regular behaviour for $\E(2)$ (left; bounded) and $\E(3)$ (right; corkscrew).}
\label{fig.E2E3}
\end{center}
\end{figure}

For chaotic shape dynamics on $X$, it is convenient to split off the compact group part of the dynamics and rewrite the system \eqref{e.iso} in the same form as in the anisotropic case \eqref{e.aniso_x}--\eqref{e.aniso_p}, namely
\begin{align*}
\dot y = f_h(y),
\quad
\dot p = \phi(y)\; ,
\end{align*}
where $y=(x,A)$, $f_h(x,A)=(f(x),Ah(x))$ and $\phi(x,A) = A v(x)$. Then
\begin{align}
p(t) = \int_0^t\phi(y(s))\, ds\; .
\label{e.py}
\end{align}
As in the anisotropic case, let $\Lambda\subset X$ be a chaotic attractor with ergodic invariant probability measure $\mu$. The natural invariant probability measure for the $y$ dynamics on $\Lambda\times \SO(d)$ is the product $m=\mu\times\nu$ where $\nu$ denotes Haar measure on the rotation group $\SO(d)$. It follows from ~\cite{FMT05} that $m$ is typically ergodic if the shape dynamics is strongly or weakly chaotic. Since $\int_{{\bf SO}(d)}A\,d\nu=0$, the time average of $\phi(x,A)=Av(x)$ is zero: 
\[
\int_{\Lambda\times{\bf SO}(d)} \phi\, dm 
= \int _{{\bf SO}(d)}A\, d\nu \int_\Lambda v\, d\mu =0\; .
\]
Hence ergodicity of $m$ implies that $\frac1t p(t)\to0$ and so there is no linear drift in the isotropic case~\cite{NicolMelbourneAshwin01,AshwinMelbourneNicol01}. 

Strong convergence to Brownian motion as in \eqref{eq-ASIP} has not been proved in the isotropic case. However for strongly chaotic dynamics, weak convergence to $d$-dimensional Brownian motion holds, and the central limit theorem follows. The covariance matrix is now diagonal due to the rotation symmetry with $\Sigma=\sigma^2I_d$ and typically $\sigma^2>0$. For rigorous results we refer the interested reader to \cite{NicolMelbourneAshwin01,FMT03,MelbourneTorok04,MelbourneNicol04}.\\

The main new result of this paper concerns a nonlinear Huygens principle for weakly chaotic systems in isotropic spatially extended media. We establish that weakly chaotic isotropic systems with odd spatial dimensions exhibit superdiffusion, whereas the superdiffusion is suppressed in systems of even dimensions. This dichotomy can be motivated from our results on regular dynamics in isotropic media as stated above. Recall that in the nonisotropic case, the anomalous diffusion is caused by the combination of laminar phases near ``sticky'' pockets of regular dynamics interspersed with intermittent chaotic bursts. However, we have seen that in even dimensions these laminar regular phases are averaged out in the isotropic case due to the rotation symmetry. Hence the mechanism for anomalous diffusion is no longer present.  We deduce that for $d$ even, weak chaos leads to Brownian behaviour just as for strong chaos. Whereas for $d$ odd, the laminar regular phase survives the effect of the rotation symmetry and we expect weak convergence to a L\'evy process.

Our predictions are supported both by the above theoretical justification and by numerical investigations described below.   (A rigorous mathematical proof is  the subject of ongoing work.)  The ingredients for the theoretical justification are summarised in a separate paragraph below. The numerical experiments also provide a useful visualisation of these phenomena. Figure~2 (bottom) presents results for an isotropic medium with $d=3$ where our theory predicts anomalous diffusion.   The computed solution behaves like a combination of Brownian motion corresponding to the intermittent chaotic bursts (as in the strongly chaotic case) together with L\'evy flights corresponding to the sticky pockets of regular dynamics. In contrast, in Figure~3 (bottom) (isotropic medium $d=2$), the computed solution behaves like a Brownian motion as before during the chaotic bursts, but the L\'evy flights are suppressed during the regular phases. The different behaviour of solutions during the regular phases, cf.\ equations~\eqref{e.peven} and~\eqref{e.podd} for $d$ even and $d$ odd respectively, is seen to be the explanation for our nonlinear Huygens principle for anomalous diffusion.    


\section{Strong and weak chaos}

In this section, we provide as promised the definition for strong and weak chaos used throughout this paper. Strongly chaotic systems include Anosov flows (Gallavotti-Cohen chaotic hypothesis \cite{GallavottiCohen95}) and uniformly hyperbolic (Axiom~A) attractors. A more general class of flows are those with a Poincar\'e map modelled by a Young tower with exponential decay of correlations~\cite{Young98}. These include H\'enon-like attractors, Lorenz-like attractors and Lorentz gas models. Even more generally, we consider situations where the Poincar\'e map is modelled by a Young tower with subexponential decay of correlations~\cite{Young99}, distinguishing between the cases where the decay rate is summable and nonsummable. For us, strongly chaotic flows are precisely those corresponding to the summable case.  (This terminology is not completely standard; many authors refer to the entire subexponential case as being weakly chaotic since Lyapunov exponents vanish.  However, as evidenced by the results described in this paper, in many respects such systems behave identically to the exponential case provided the decay is summable, and it is the boundary between summable and nonsummable that is significant.)

Roughly speaking, weakly chaotic flows are those corresponding to the nonsummable case, but there is an extra requirement that decay rates are regularly varying functions\footnote{A function $\ell(x)$ is {\em slowly varying} if $\ell(\lambda x)/\ell(x)\to1$ as $x\to\infty$ for all $\lambda>0$; examples are functions that are asymptotically constant and powers and iterates of logarithms. A function of the form $\ell(x)x^q$ is {\em regularly varying}.}. This is not simply a technical hypothesis; regular variation of tails is a necessary condition for convergence to a stable law or L\'evy process.

Note that the definition of strong/weak chaos makes assumptions on the decay of correlations for the Poincar\'e map, but not for the attractor $\Lambda$ itself. This is important since even Anosov flows are not necessarily mixing and there are mixing uniformly hyperbolic flows with arbitrarily slow decay of correlations~\cite{Ruelle83,Pollicott85}. In particular, mixing properties for the Poincar\'e map do not necessarily pass to the flow. In contrast, convergence to a Brownian motion or a L\'evy process does pass to the flow~\cite{MelbourneTorok04,Gouezel07,MelbourneZweimueller13}.
 

\begin{figure}[t]
\begin{center}
\centerline{
\includegraphics[width=.2\textwidth]{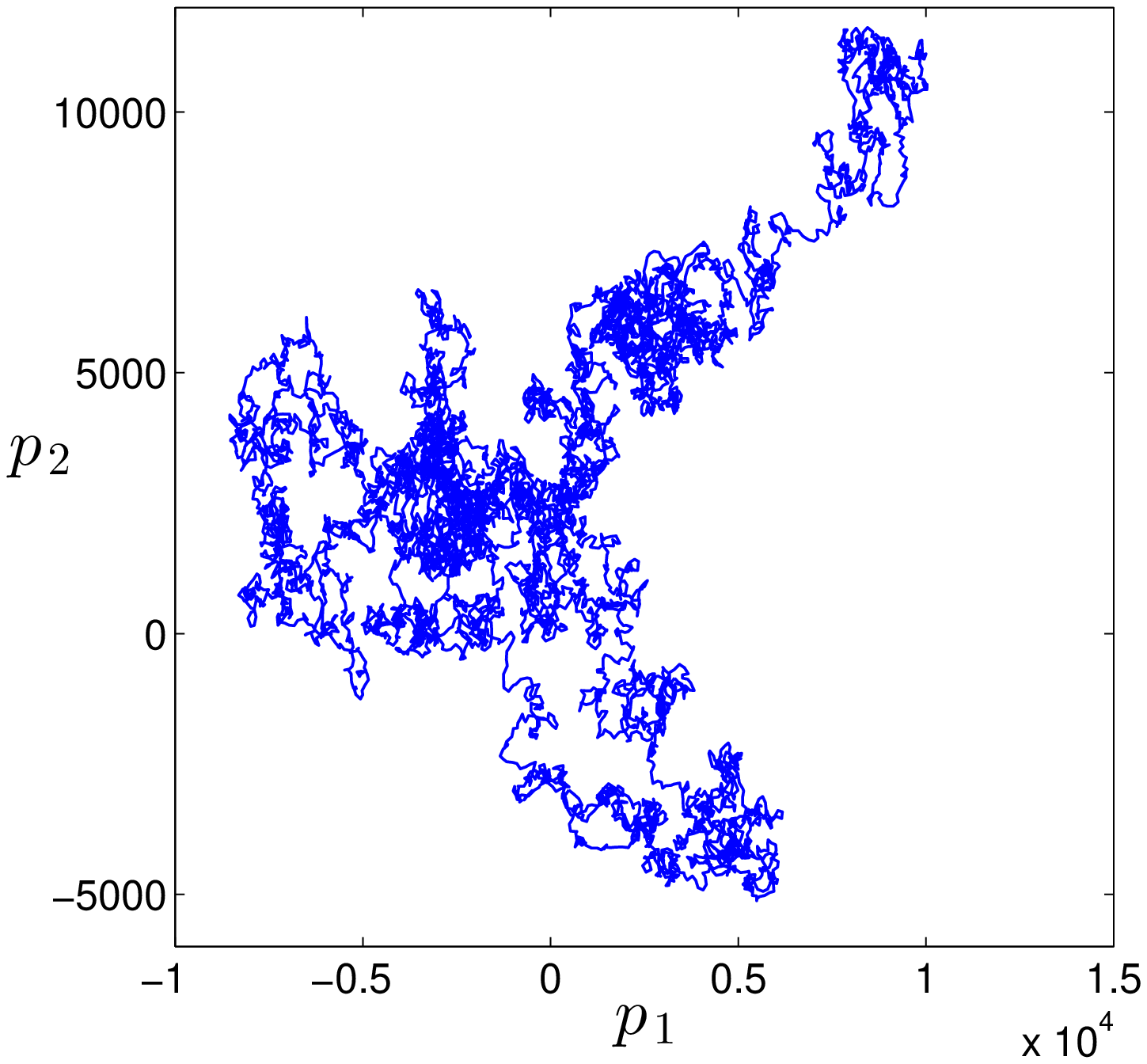}
\includegraphics[width=.2\textwidth]{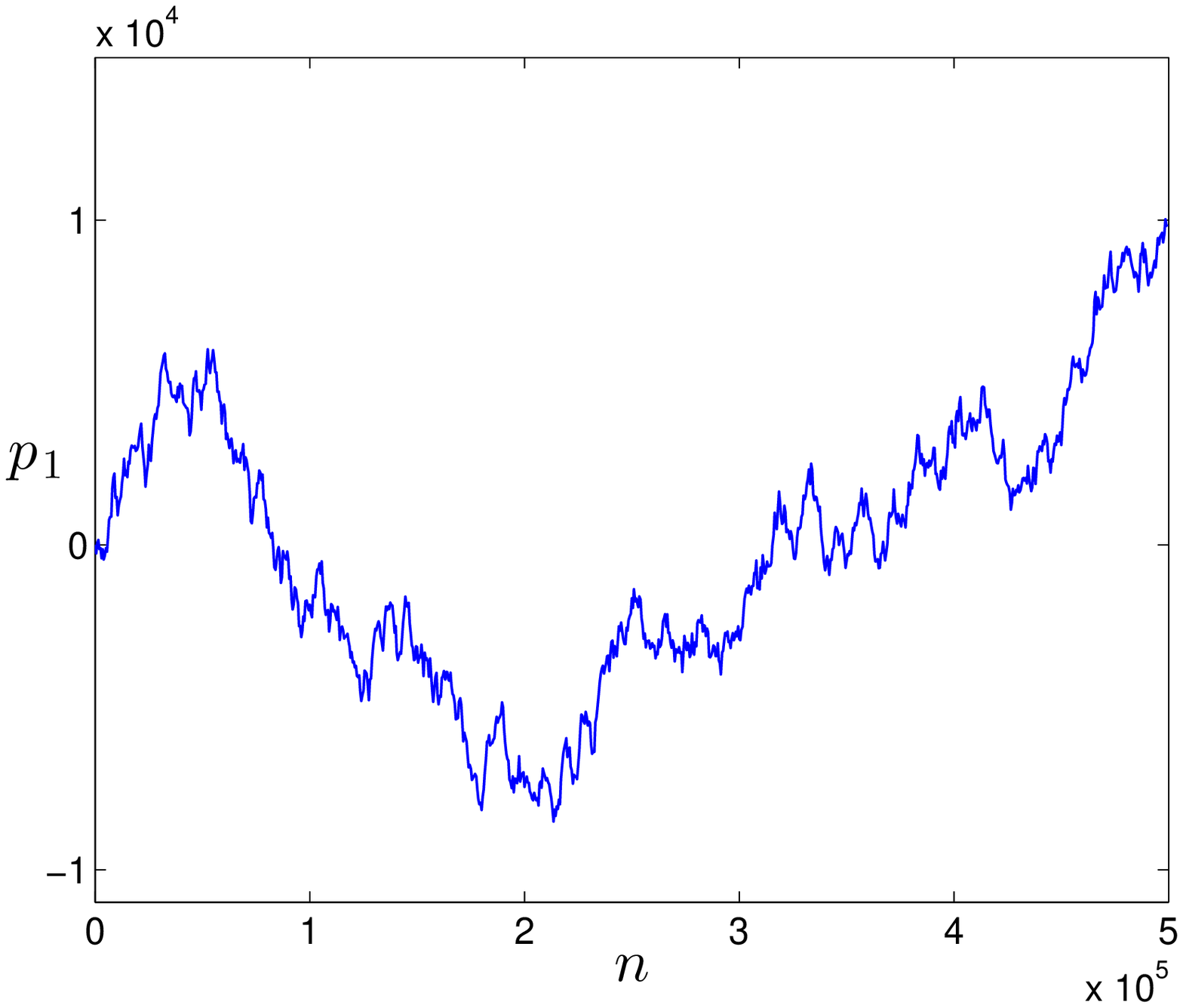}
}
\centerline{
\includegraphics[width=.2\textwidth]{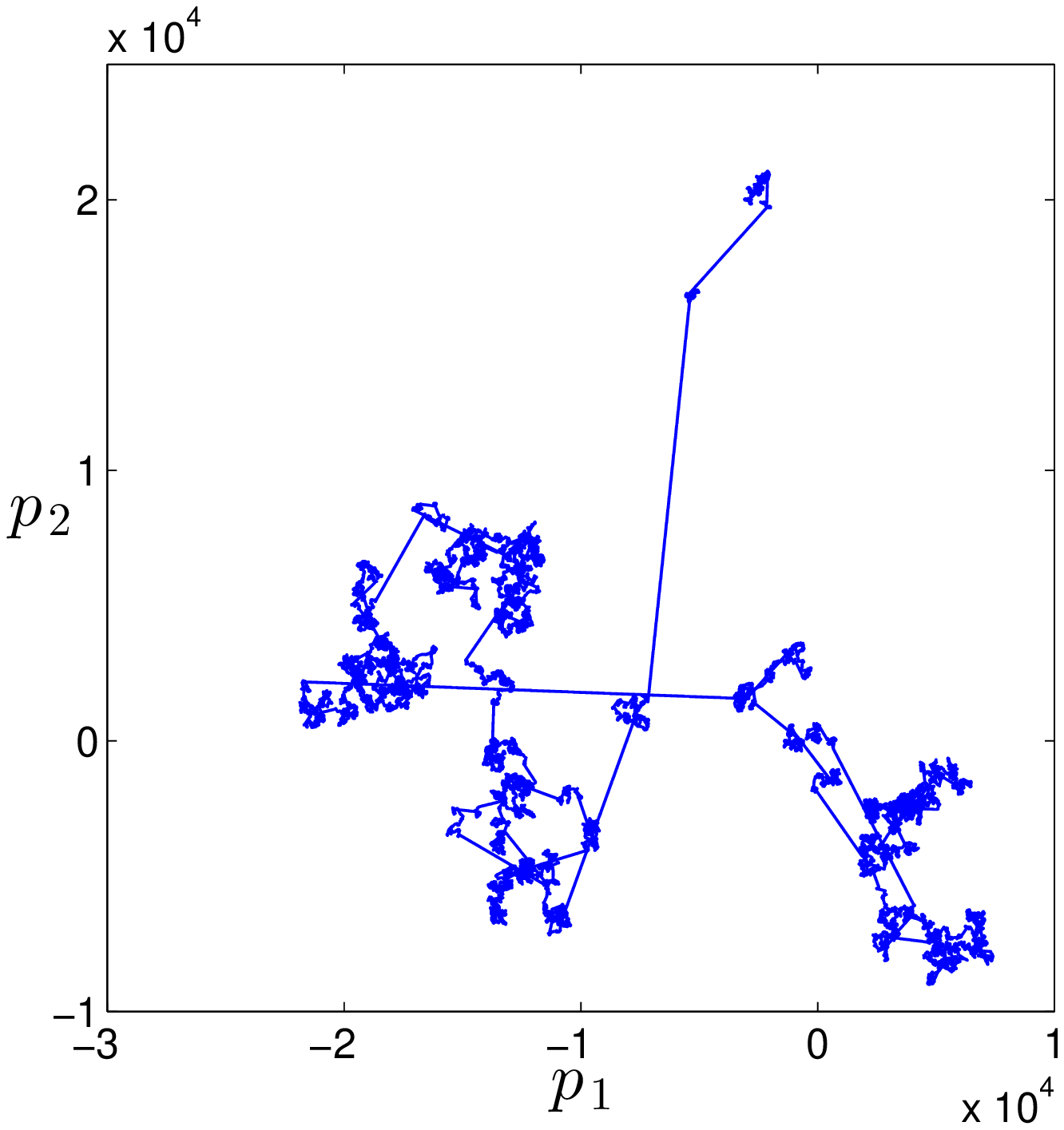}
\includegraphics[width=.2\textwidth]{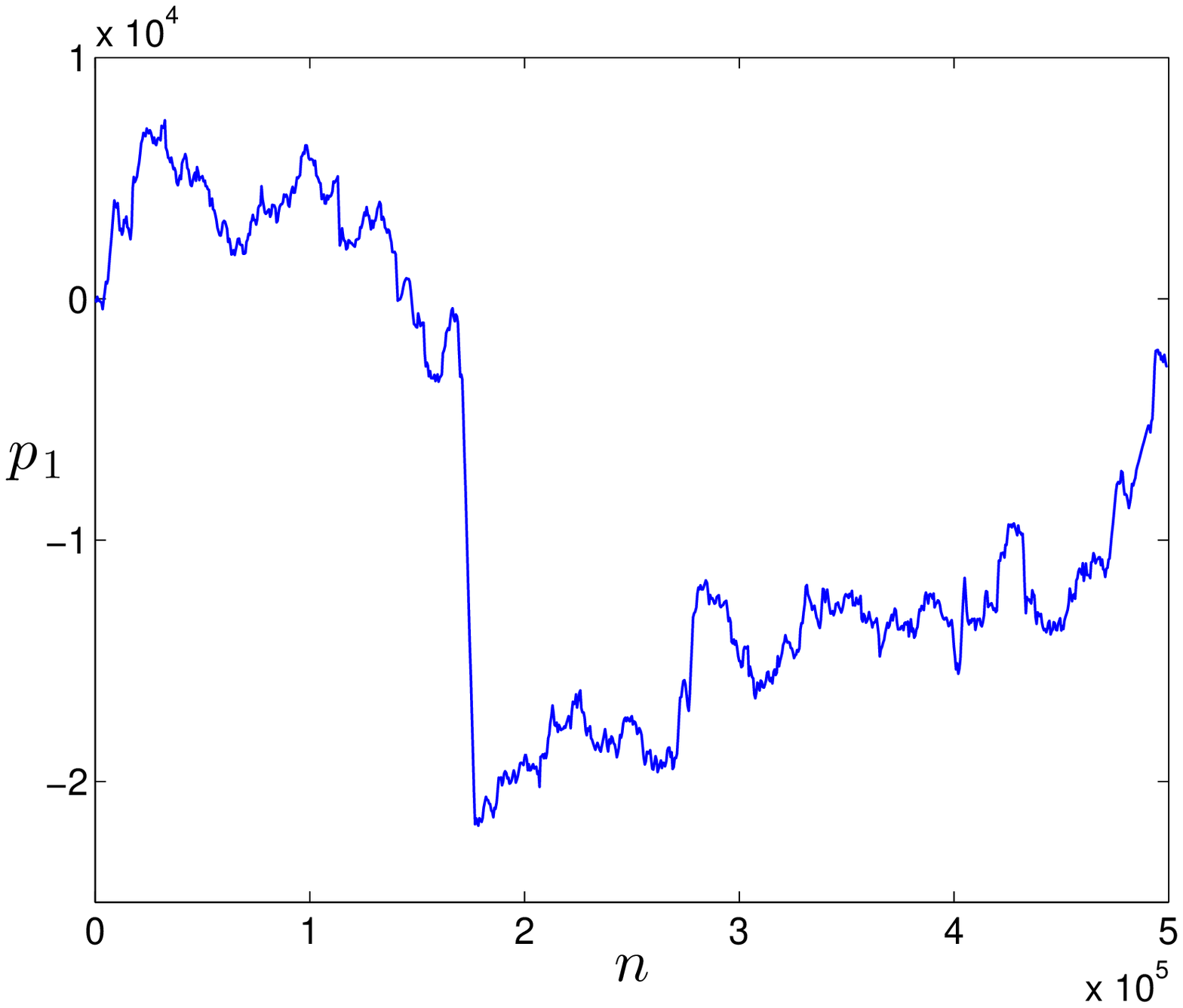}
}
\caption{Isotropic case in odd dimensions: Translation variables $p_1$ and $p_2$ of an $\E(3)$ skew product driven by the Pomeau-Manneville map \eqref{e.PM}. Top: Strongly chaotic case with $\gamma=0.2$: Brownian motion. Bottom: Weakly chaotic case with $\gamma=0.7$: L\'evy process.}
\label{fig.E3extPQ}
\end{center}
\end{figure}

\begin{figure}[t]
\begin{center}
\centerline{
\includegraphics[width=.2\textwidth]{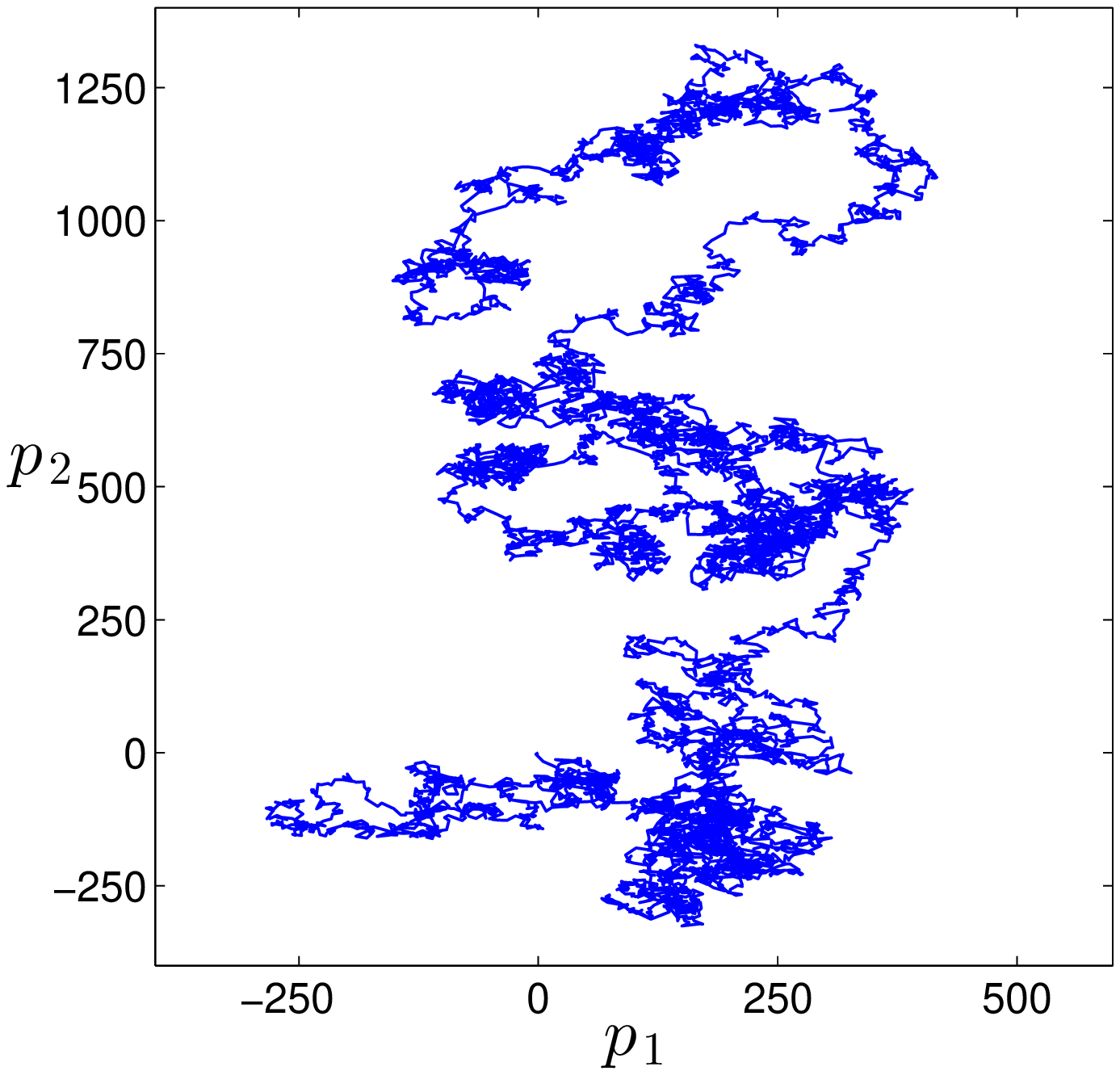}
\includegraphics[width=.2\textwidth]{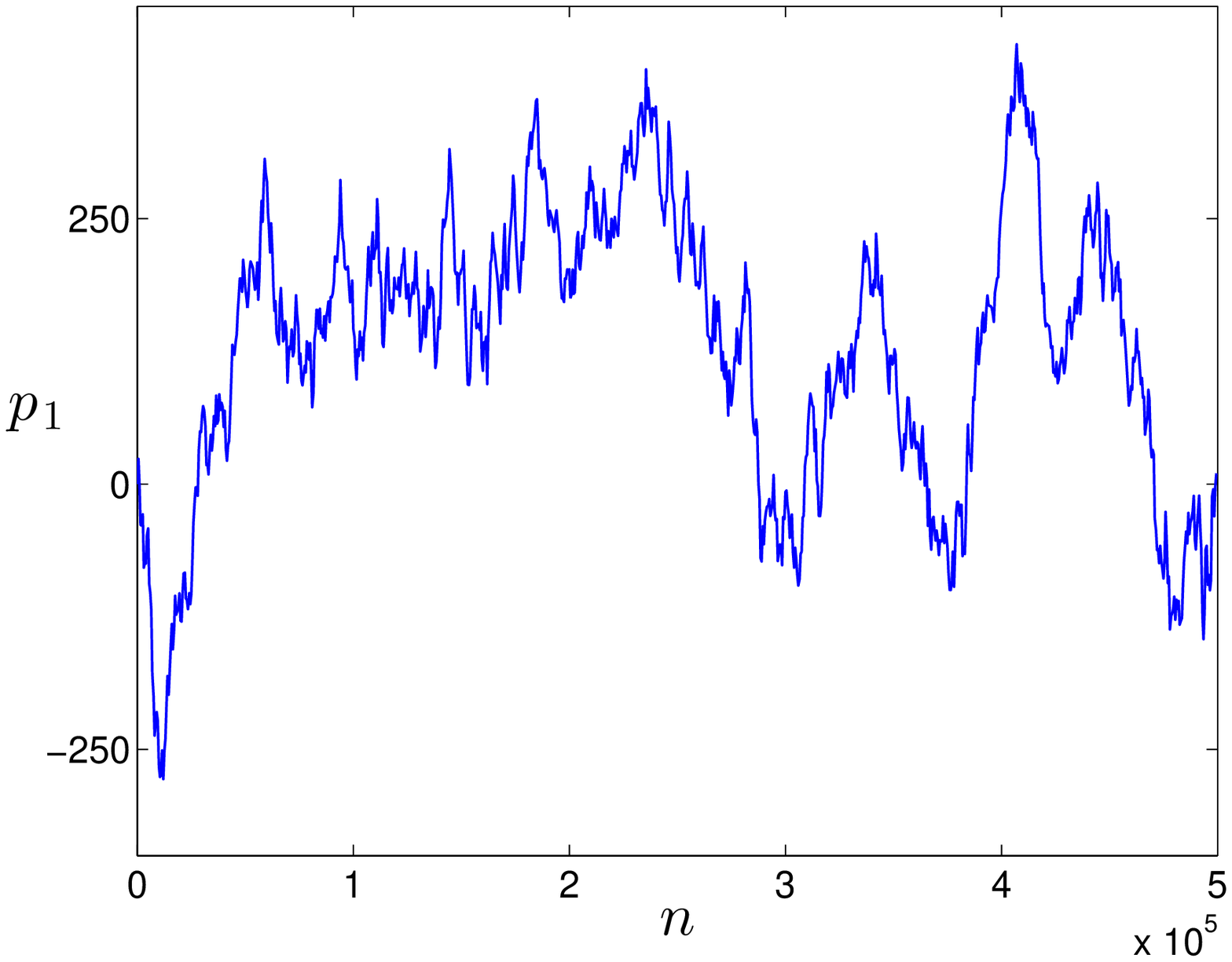}
}
\centerline{
\includegraphics[width=.2\textwidth]{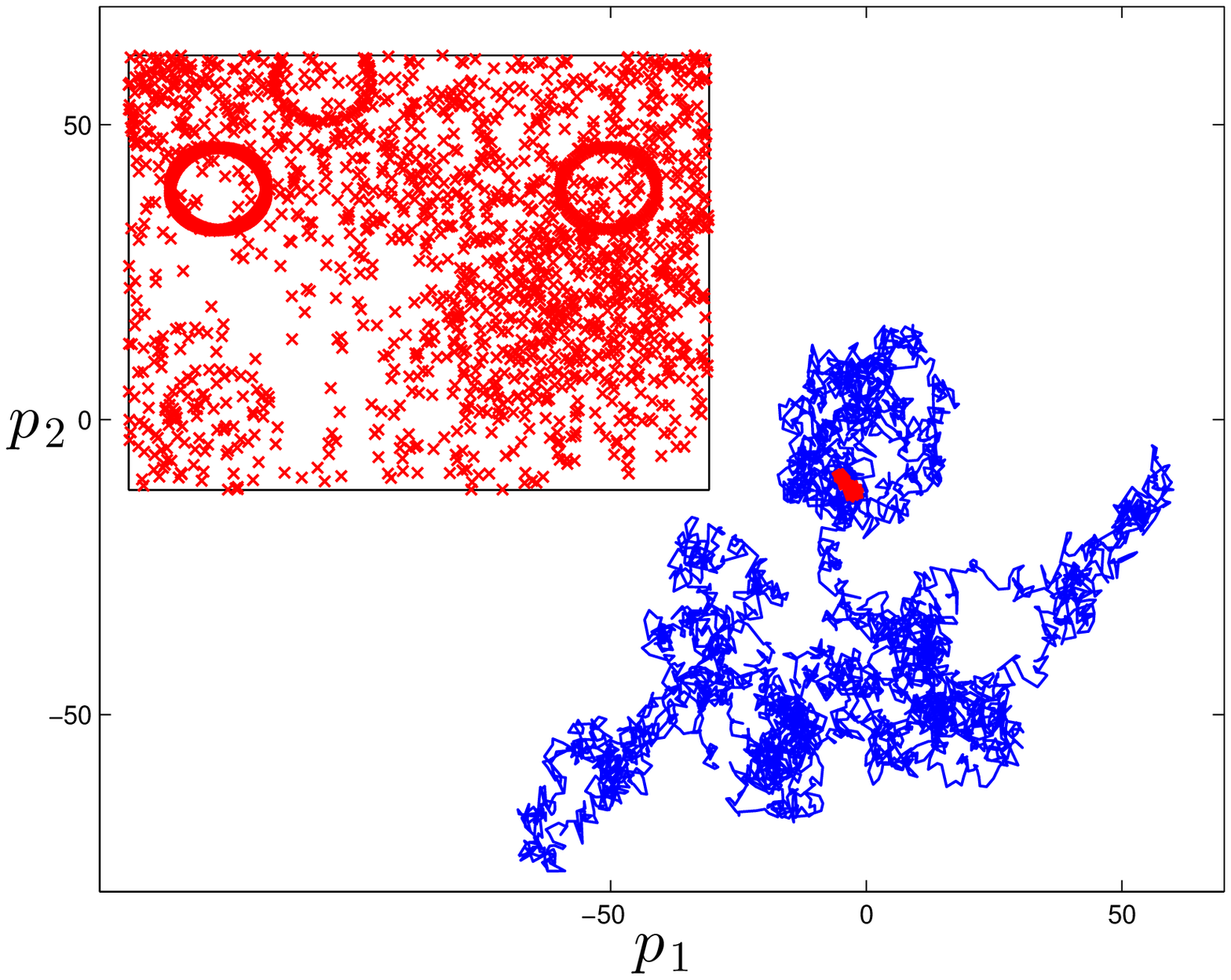}
\includegraphics[width=.2\textwidth]{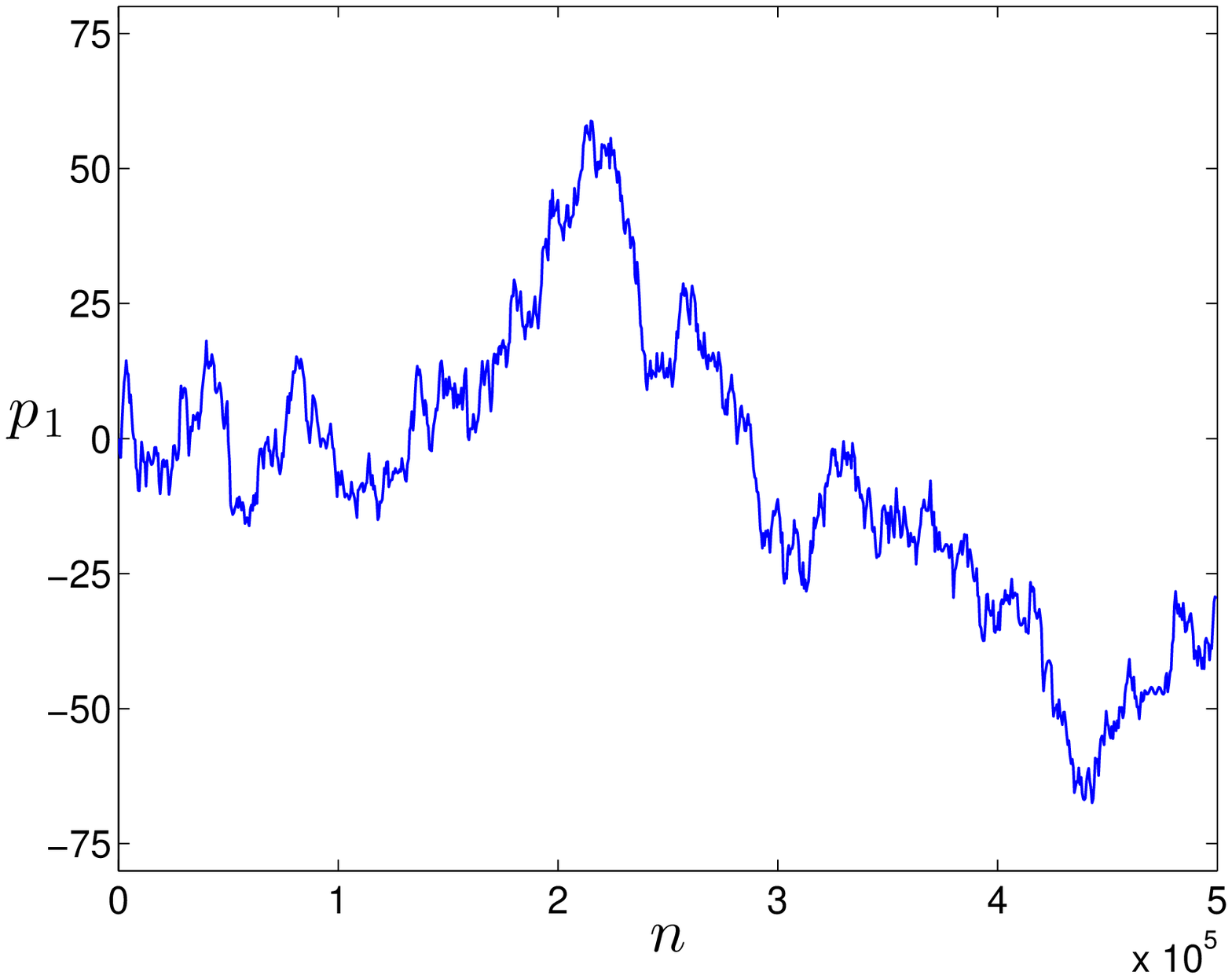}
}
\caption{Isotropic case in even dimensions: Translation variables $p=(p_1,p_2)$ of an $\E(2)$ skew product driven by the Pomeau-Manneville map \eqref{e.PM} exhibiting Brownian motion. Top: Strongly chaotic case with $\gamma=0.2$. Bottom: Weakly chaotic case with $\gamma=0.7$ with inset zooming into a smaller area illustrating the mechanism of suppression of anomalous diffusion in even space dimensions.}
\label{fig.E2extPQ}
\end{center}
\end{figure}

\section{Theoretical Justification}
The individual ingredients comprising the theoretical justification for the results in this paper are standard in various subsections of the scientific community, but may seem nonstandard when taken collectively.  Hence it seems worth summarising the ingredients here.
\begin{itemize}
\item   The underlying pattern-forming system (which could for example be a partial differential equation or the physical system itself) is decomposed in the skew product form~\eqref{eq-skew} (see for example \cite{GolubitskyStewart}).   
\item  The assumed form of the shape dynamics $x(t)$ leads to  various types of equation $\dot g(t)=g\xi(x(t))$ for the group dynamics.  There is  a well-established theory when the shape dynamics is steady or periodic~\cite{Field80,Krupa90,AshwinMelbourne97}, and there are numerous results also in the case when the shape dynamics is chaotic~\cite{NicolMelbourneAshwin01,FMT03,MelbourneNicol04},
\item 
The ``regular dynamics'' entries in Table~1 are immediate from~\cite{AshwinMelbourne97}. This already establishes the absence of L\'evy flights during regular 
phases in weakly chaotic systems for even-dimensional isotropic media.
\item  Passing to a Poincar\'e cross-section reduces to the conceptually and technically simpler situation of discrete time.   Results of~\cite{MelbourneTorok04} guarantee that statistical limit laws for the Poincar\'e map extend back to the continuous time setting.
\item   Young towers with exponential decay of correlations~\cite{Young98} have been shown to include most of the classical examples of strongly chaotic dynamical systems (eg H\'enon, Lorenz) extending far beyond the Anosov/Axiom A setting.  A broader class is Young towers with summable decay of correlations~\cite{Young99}.  The statistical properties of such systems have been the subject of much recent mathematical investigation, and the ``strongly chaotic'' entries in Table~1 are special cases of these rigorous results \cite{NicolMelbourneAshwin01,FMT03,MelbourneTorok04,MelbourneNicol04}.
\item  Even more recently, starting with~\cite{Gouezel04}, the statistical properties of Young towers with nonsummable tails have come under the spotlight.  The prototypical example is provided by Pomeau-Manneville intermittency maps.  Here anomalous diffusion is anticipated, and convergence to a L\'evy process is rigorously proved in~\cite{MelbourneZweimueller13}.   The bottom left entry in Table~1 is a consequence of this.
\item  The remaining entries in the bottom row of Table~1 remain conjectural from the point of view of rigorous mathematics but the theoretical justification is as follows. Intermittent dynamics is a mixture of regular phases and chaotic bursts. Based on~\cite{GaspardWang88,Gouezel04,MelbourneZweimueller13} one is led to anticipate superdiffusive behaviour with L\'evy flights corresponding to the regular phases. But we already saw that regular dynamics in isotropic media varies significantly in even and odd dimensions and that the mechanism for L\'evy flights exists only in odd dimensions.  Consequently we predict that anomalous diffusion is suppressed in even-dimensional isotropic media and exists only in odd dimensions, and this prediction is supported by the numerical experiments below. (A recent rigorous result of~\cite{PeligradWu10} shows that anomalous diffusion is indeed suppressed in the case of two dimensions under the assumption that the rotation component $h\in\SO(2)$ is constant.  The assumption on $h$ enables the use of Fourier analysis and simplifies matters significantly.  However, the heuristics behind our results do not rely on this assumption.)
\end{itemize}

In general, there is no convenient way to explicitly determine the skew-product equations~\eqref{eq-skew} from the underlying pattern-forming system.   To circumvent this, it has become standard in the Physics literature to consider lattice models for diffusion, where the underlying system is posed in physical space $\R^d$
(for example to model phase dynamics of Josephson junctions and charge-density waves
\cite{GeiselN82, GNZ85} and advection-diffusion of passive tracers in fluid flows~\cite{Pierrehumbert00}).
The Euclidean group 
is replaced by a discrete group $G$ of rotations and translations, and $X$ is identified with a fundamental domain for the action of $G$
on $\R^d$.   In this setting, the transition between the underlying equations and the skew-product system is completely transparent and explicit.   Our results hold equally in this situation, see the
Supplementary Information for further details and numerical results.


\section{Numerical results}

Intermittent dynamics is often modelled by the prototypical family of Pomeau-Manneville intermittency maps $x_{n+1}=f(x_n)$ with $f:[0,1]\to[0,1]$ given by 
\begin{align}
f(x)=
\begin{cases}
x(1+2^\gamma x^\gamma) & 0\le x\le \frac12 \\ 
2x-1 & \frac12\le x\le1 
\end{cases}
\label{e.PM}
\end{align}
where $\gamma$ is a parameter \cite{PomeauManneville80,LiveraniEtAl99}. Pomeau-Manneville intermittency maps are the prototype for the study of intermittency in finite and infinite dimensional systems where they have been used to model Poincar\'e maps (see for example \cite{BergePomeauVidal,CovasTavakol99}). If $\gamma\in[0,1)$ there exists a unique absolutely continuous invariant probability measure (SRB measure) $\mu$. When $\gamma=0$ this is the doubling map with exponential decay of correlations. For $\gamma \in (0,1)$, it is known~\cite{Hu04} that the decay of correlations is polynomial with rate $1/n^{(1/\gamma)-1}$ which is summable for $\gamma<\frac12$ and nonsummable for $\gamma \in [\frac{1}{2},1)$; according to our definition the Pomeau-Manneville map \eqref{e.PM} is strongly chaotic for $\gamma\in[0,\frac12)$ and weakly chaotic for $\gamma \in [\frac{1}{2},1)$. Note that for $\gamma>0$ the fixed point at $0$ is indifferent ($f'(0)=1$) and plays the role of the ``sticky'' regular dynamics. For $\gamma\ge\frac12$, the stickiness is strong enough to support superdiffusive phenomena. (In the borderline case $\gamma=\frac12$, there is still weak convergence to Brownian motion but with anomalous diffusion rate $\sqrt{t\log t}$ in the anisotropic case~\cite{Gouezel04,DolgopyatEtAl08} and similarly for isotropic systems in odd dimensions.). 
We now present numerical results for skew products representing the anisotropic and isotropic cases (both even and odd dimensional) in the strongly chaotic and weakly chaotic regimes. We make the obvious modifications to the continuous time case described earlier so the skew product is now discrete in time.

\noindent
{\bf Anisotropic case}. 
Here $G$ is the translation group $\R^d$, and the discrete time skew product for the Pomeau-Manneville map \eqref{e.PM} reads as 
$(x,p)\mapsto (f(x), p+\phi(x))$ so that
\begin{align}
p(n) = \sum_{j=0}^{n-1} \phi(x_j)\;.
\label{e.PMR}
\end{align} 
In Figure~4, we take $d=1$ and $\phi(x)=1+x$. The translation coordinate $p(n)$ exhibits linear drift $cn$, where $c=\int_{[0,1]}\phi\,d\mu$, for both strongly and weakly chaotic dynamics. Passing into the comoving frame, the distinction between Gaussian and L\'evy type fluctuations becomes apparent. 

Since the fixed point at $0$ is indifferent for $\gamma>0$, an initial condition $x_0$ that starts close to $0$ remains close to $0$ for a large number of iterates, say $N_0$, so that $\phi(x_j)$ is roughly of size $\phi(0)$ for $j=0,1,\ldots,N_0$. Hence in the comoving frame $p(n)$ exhibits approximately linear growth, $p(n)\approx (\phi(0)-c)n$, for $n\le N_0$. In particular, the small jumps $\phi(0)-c$ accumulate into a large jump. This is akin to a particle's ballistic motion with constant velocity $\phi(0)-c$, a common picture of anomalous diffusion \cite{SolomonEtAl1994}.

In the strongly chaotic case $\gamma\in[0,\frac12)$, these large jumps are too rare to cause anomalous diffusion and, as explained above, it follows from~\cite{MelbourneNicol09} that 
\[
p(n)=cn+W(n)+O(n^{\frac12-\lambda})\qquad {\rm{a.e.}}
\]
(cf.\ equation~\eqref{eq-ASIP}) so the dynamics in the comoving frame is Brownian-like. However, in the weakly chaotic case $\gamma\in(\frac12,1)$, Gou\"ezel~\cite{Gouezel04} demonstrated that the large jumps correspond to L\'evy flights. Once the trajectory $x_n$ moves away from $0$, the values of $\phi$ fluctuate erratically yielding Brownian-like motion for $p(n)$. These two effects combine~\cite{Gouezel04,MelbourneZweimueller13} to
produce a L\'evy process with diffusion rate $t^\gamma$. 

\noindent
{\bf Isotropic case}. We illustrate the nonlinear Huygens principle by which weakly chaotic dynamics causes anomalous diffusion in isotropic media with odd space dimensions but normal diffusion in even space dimensions. The $\E(2)$ skew product with rotations $\theta$ and translations $p=(p_1,p_2)$ can be represented as
\[
(x,\theta,p(x))\mapsto (f(x),\theta+h(x),p+e^{i\theta}v(x)),
\]
where $h$ and $v$ takes values in $\R$ and $\C$ respectively. Hence
\begin{align}
p(n) = \sum_{j=0}^{n-1} e^{i\theta_j} v(x_j) 
\quad  
\theta_j = \sum_{k=0}^{j-1} h(x_k).
\label{e.PME2}
\end{align}
For the numerics, we choose $v(x)=1+x$ and $h(x)=c_0\neq0$. In Figure~3 we show plots of the translation variables in the $(p_1,p_2)$-plane. We also show the process $p_1(n)$ as a function of time. The diffusive behaviour is seen to be normal in both the strongly and weakly chaotic cases. The mechanism by which anomalous diffusion is suppressed in even dimensions is nicely illustrated: the long laminar phases which in the anisotropic case give rise to large excursions are now bounded in both the strongly and weakly chaotic case. This is seen in the intermittent circular motion in the inset of Figure~3. When combined with the Brownian like behaviour caused by the chaotic bursts, this leads to an overall Brownian behaviour. It is clear from \eqref{e.PME2} that if we were to remove the rotation near the fixed point by setting $h(x)\equiv0$ for $x$ near zero, then the translation variables would exhibit the same behaviour as in the anisotropic case (cf.~\eqref{e.PMR}). In particular, this would yield a L\'evy process in the weakly chaotic case.\\ The corresponding plots for the $\E(3)$ skew product are shown in Figure~2.  The anomalous L\'evy type behaviour is clearly visible in the weakly chaotic case.  See SI for details on the skew product equations used for the numerics.

\begin{figure}
\begin{center}
\centerline{
\includegraphics[width=.25\textwidth]{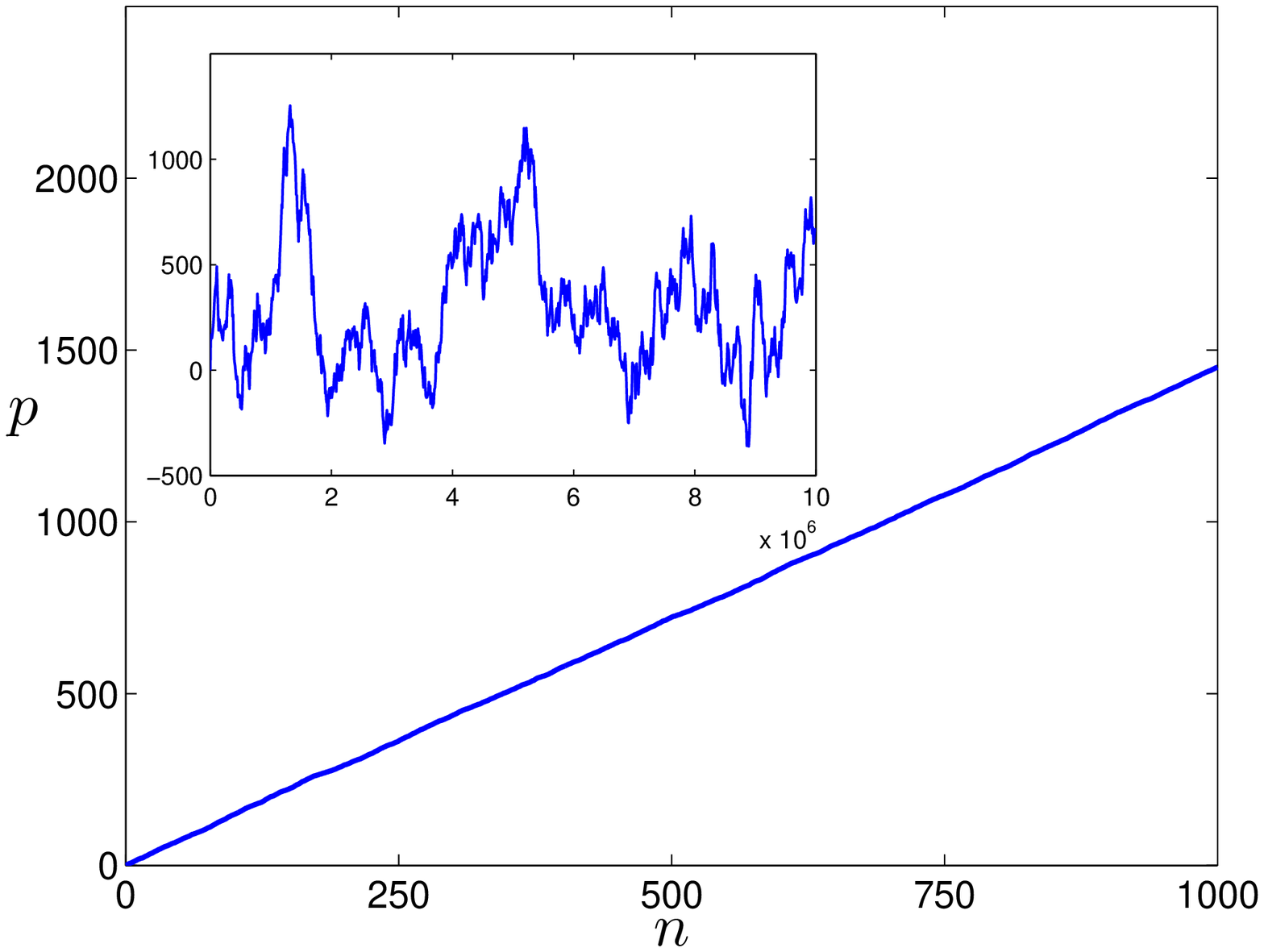}
\includegraphics[width=.25\textwidth]{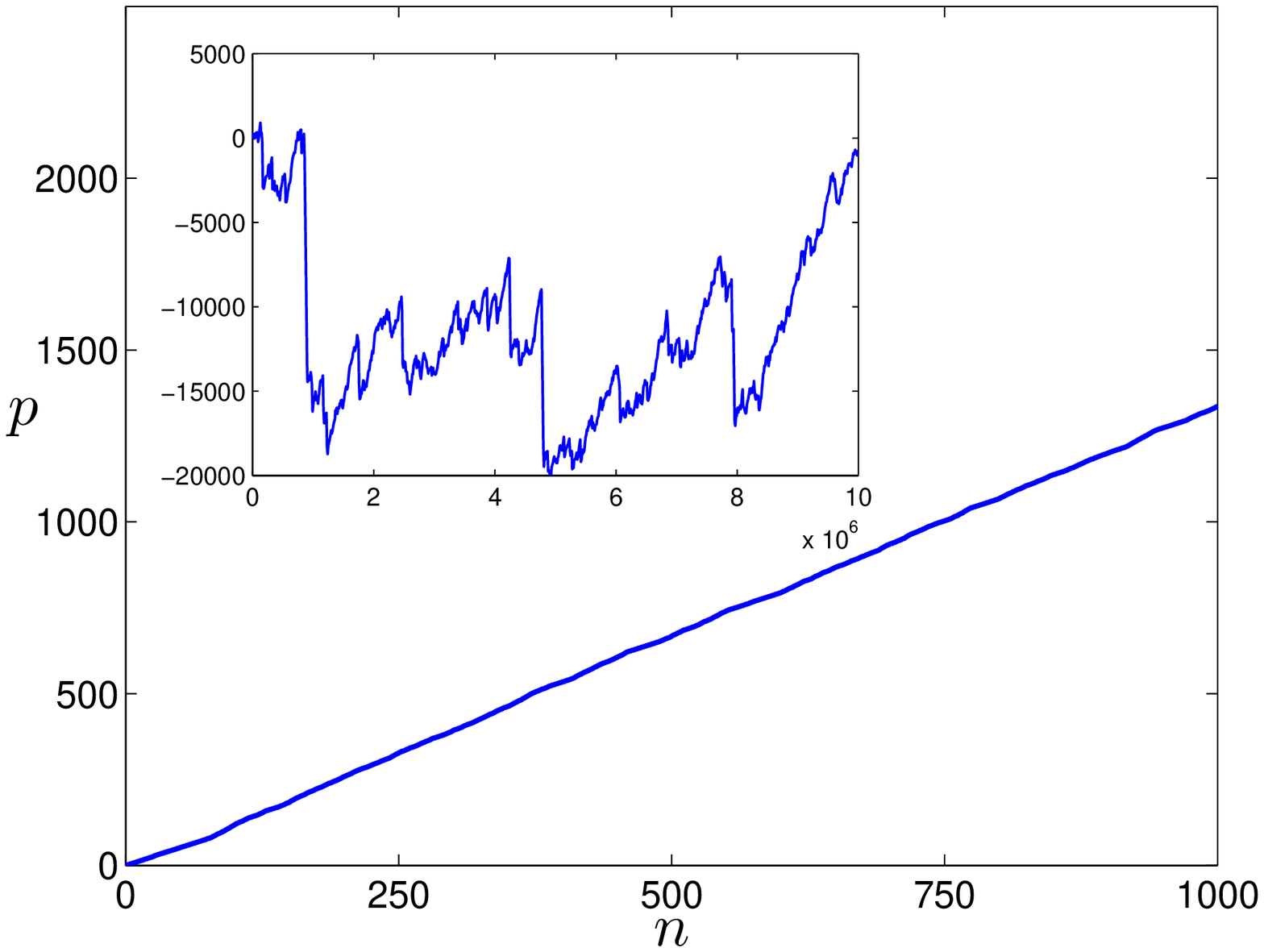}}
\caption{Anisotropic case: Translation variable $p$ as a function of time for an $\R$ skew product driven by the Pomeau-Manneville map \eqref{e.PM}. Left: Strongly chaotic case with $\gamma=0.2$: linear drift with superimposed Brownian motion. Right: Weakly chaotic case with $\gamma=0.7$: linear drift with superimposed L\'evy process. To highlight the diffusive behaviour, we eliminate the linear drift by subtracting the mean from the data. Note that in the weakly chaotic case the L\'evy process is asymmetric; the flights are downwards since $\phi(0)-c$ is negative.}
\label{fig.RextP}
\end{center}
\end{figure}


\section{Summary and discussion}

We provided a universal view on the type of diffusive behaviour which can be expected in spatially extended systems with symmetry. In doing so we proposed a new definition of weak chaos as the boundary between summable and non-summable correlations. In contrast to the previous view of the onset of weak chaos as the demarcation line of exponential and sub-exponential decay of correlations, our definition allows for a distinction between normal and anomalous behaviour. Using this definition we contributed to the understanding of diffusive behaviour in isotropic media establishing a nonlinear Huygens principle whereby superdiffusion occurs naturally in odd dimensions but not in even dimensions.

The phenomenon of anomalous diffusion has attracted a lot of interest in the last decade, with applications ranging from  the motion of metal clusters and large molecules across crystalline surfaces \cite{SanchoEtAl04}, conformational changes in proteins \cite{ReuveniEtAl10}, migration of epithelial cells \cite{DieterichEtAl08}, diffusion in plasma membranes of living cells \cite{WeigelEtAl11}, finance ratios \cite{Podobnik08} to foraging strategy of animals \cite{Fernandez03,GetzSaltz08,BartumeusLevin08} to name but a few. The mechanisms for anomalous diffusion in these papers are model-dependent, specific to the details of the geometry of the various situations. In contrast, ours is a universal perspective (analogous to the classical Huygens principle) driven only by the ambient symmetry and the degree of chaoticity of the underlying dynamics. Our theory sets out the general conditions under which anomalous diffusion can be expected in spatially extended systems with symmetry and may be viewed as a prediction and interpretation of superdiffusive behaviour in future experiments. \\

It is of interest to consider diffusive and superdiffusive behaviour in systems with different kinds of spatial symmetry.    For systems with rotation $\SO(2)$ symmetry, for example rotating convection in the plane or on a sphere, trajectories lift naturally to the universal covering group $\R$, and we obtain the same results as for the anisotropic case with $d=1$.     

For systems with $\O(2)$ symmetry (rotations and reflections), it is interesting to consider the effect of the reflection symmetry. There are two entirely different scenarios: (i) There are two disjoint attractors interchanged by reflections and each behaves as in systems with $\SO(2)$ symmetry. (ii) The attractor is invariant (setwise) under reflections and the linear drift vanishes; we predict bounded trajectories for regular dynamics and Brownian motion for strongly chaotic dynamics. But in the case of weakly chaotic dynamics, we expect diffusion or superdiffusion depending on whether the sticky regular phase is invariant or not under reflections. Further, in scenario~(ii) the Brownian motions and L\'evy processes are symmetric.

Finally, we mention an open problem about systems on a sphere with $\SO(3)$ rotation symmetry (again these could be reaction diffusion equations or convection problems). This time, there is no elementary method for passing to a noncompact group where it makes sense to speak of unbounded growth of trajectories. Our expectation is that locally the results are similar to those in the unbounded plane ($\E(2)$ symmetry) but it is unclear how to make such a statement precise.


\begin{acknowledgments}
GAG acknowledges funding from the Australian Research Council. The research of IM was supported in part by EPSRC Grant EP/F031807/1. IM is grateful for the hospitality of the University of Sydney where most of this research was carried out.

\end{acknowledgments}


\end{article}


\end{document}